# A note on comparison of scientific impact expressed by the number of citations in different fields of science

**Igor Podlubny**

*Department of Applied Informatics and Process Control*
*Faculty of B.E.R.G., Technical University of Kosice*
*B. Nemcovej 3, 04200 Kosice, Slovak Republic*
*Phone: +421 55 6025179, Fax: +421 55 6025190, +421 55 6336618*
E-mail: igor.podlubny@tuke.sk

**Abstract:** Citation distributions for 1992, 1994, 1996, 1997, 1999, and 2001, which were published in the 2004 report of the National Science Foundation, USA, are analyzed. It is shown that the ratio of the total number of citations of any two broad fields of science remains close to constant over the analyzed years. Based on this observation, normalization of total numbers of citations with respect to the number of citations in mathematics is suggested as a tool for comparing scientific impact expressed by the number of citations in different fields of science.

## 1 Introduction

Number of citations is usually considered as one of important indicators of the scientific impact of a scientist in his/her particular field. This criterion can be easily used in each particular field, when two mathematicians (or two physicists, or two chemists, or two medical researchers, or two engineers, etc.) are compared. This comparison is used by the Thomson ISI (the Institute for Scientific Information) for compiling various lists, like ISIHighlyCited.Com [2], arranged by scientific field.

A more difficult problem arises when we have to compare two scientists working in *different* fields, for example, a mathematician and a chemist. The difficulty is underlined by the fact that even the most prolific author of citation analysis, Dr. E. Garfield, used only absolute figures for compiling lists of scientists with the highest impact – see, for example, the list in [6], where we cannot see any mathematician, engineer, or a specialist in social sciences. The same approach (total numbers of citations) is used also in [7], where one can observe the same absence of mathematicians, physicists, engineers, etc.

It is obvious that we cannot compare total numbers of citations – it is well known that in absolute figures there is much less citations in mathematics than in chemistry, but a mathematician with a relatively low total number of citations can have higher impact in mathematics than a chemist with a larger number of citations in chemistry. The question, therefore, is: *is it possible to compare two scientists working in different fields of science on the basis of their citation numbers?* Surprisingly, the author of this article could not find any answer to this seemingly natural question in the available literature.

The answer suggested in this article is: *yes, it is possible, with the help of a certain normalization of their respective numbers of citations.* The proposed approach is described below.

## 2 The Data

In a recent publication of the National Science Foundation the distribution of scientific citations of the U.S. scientific and engineering articles across wide fields of science in 1992, 1994, 1996, 1997, 1999, and 2001 was published (see [1], Chapter 5, Table 5-27 on page 5-50). The sources for the data appearing in that table were the Science Citation Index (SCI) and the Social Sciences Citation Index (SSCI).

## 3 The Law of the Constant Ratio

The data in the NSF table for the distribution of scientific citations led me to the observation that *the ratio of the number of citations in any two fields of science remains close to constant*.

For example, for clinical medicine and physics we have the ratio close to 4:

(1992) 475793 / 137922 = 3.44972521

(1994) 516665 / 141653 = 3.64739893

(1996) 554332 / 138417 = 4.0047971

(1997) 574859 / 131958 = 4.35637854

(1999) 584330 / 125968 = 4.63871777

(2001) 589762 / 120593 = 4.89051603

Similarly, for engineering and mathematics we obtain the ratio close to 5:

(1992) 32680 / 6858 = 4.76523768

(1994) 35189 / 6631 = 5.30674106

(1996) 33664 / 6961 = 4.83608677

(1997) 32958 / 6418 = 5.13524462

(1999) 34001 / 7520 = 4.52140957

(2001) 36809 / 7794 = 4.72273544

The same observation holds for any pair of fields of science in the Table 5-27 of the NSF 2004 report.

It is worth noting that a similar law of the constant ratio in citation analysis is known for the number of publications and the number of citations processed by ISI – it gives the so-called Garfield's constant [5].

## 4 Normalization

Based on the observed law of the constant ratio, we can normalize all scientific fields by computing the ratio of the number of citations in each field to the number of citations in mathematics (the smallest number of citations among all fields). The results are shown in Table 1 in the columns titled *"ratio to maths"*; numbers are rounded to integers. In such a form the law of the constant ratio is even more obvious.

The average ratio of citation number to the number of citations in mathematics is given in a dedicated column in Table 1.

# 5 Comparing different fields of science

Using the suggested normalization of the citation data provided by the SCI (Science Citation Index), we could – to some extent – compare the relative scientific impact of research institutions and maybe even individual scientists working in different fields of science.

**Example 1.**

- Q: Who has higher impact in his field: a physicist with 70 citations or an engineer with 20 citations?

- A: In normalized units, the physicist's impact is 70:19=3.68, while the engineer's impact is 20:5=4. Therefore, the engineer has slightly higher impact in his field than the physicist in his one (although it is not clear at all from their total numbers of citations).

**Example 2.**

- Q: How many citations can be considered as equivalent for mathematics, chemistry, physics, and clinical medicine?

- A: According to Table 1, one citation in mathematics roughly corresponds to 15 citations in chemistry, 19 citations in physics, and 78 citations in clinical medicine. In other words, 250 citations in mathematics can be considered as roughly equivalent to 3750 citations in chemistry, 4750 citations in physics, and 19500 citations in clinical medicine.

# 6 Conclusion

In conclusion, the following could be mentioned.

For the proper interpretation of the above observation, it may be important that the total number of citations in the analyzed data is almost stabilized. The additional analysis of the impact of the growing number of publications and citations from emerging regions like Eastern Europe or Asia on the constants listed in Table 1 needs further data, which should include those regions.

It seems that the *law of the constant ratio*, described in this brief note, gives reasonable results and can be used in average as a macro criterion for comparing the scientific impact of research institutions belonging to different fields of science. At micro level, it can be used (with care) for comparing scientists with low or average scientific impact from the viewpoint of citations of their works. Even in this case it is necessary to take into account that the distribution of citations in small sets of documents is usually irregular. In case of large numbers of citations it will probably need some correction, since the ratios of peaks in different fields of science do not necessary copy the ratios shown in Table 1.

It may seem that by considering only total numbers of citations in various fields of science we do not take into account the fact that the numbers of scientists working in those fields also differs significantly, as well as the total number of publications in those fields. However, the reality is the opposite: the smaller number of citations, for example, in mathematics compared to biomedicine, simply reflects the fact that the number of articles in mathematics is also smaller than the number of articles in biomedicine, that there is less people publishing in mathematics than in biomedicine, and that the average length of the reference list in mathematics is less than the average length of the reference list in biomedicine. Therefore, the differences in the number of people and in the number of publications in different fields are taken into account implicitly through the total numbers of citations produced by those people in those publications.

The approach suggested in this article can bring mathematicians, engineers, and other "less visible" scientists to the multidisciplinary lists of high-impact scientists, thus correcting the approach used in [6] and other similar lists [7].

Finally, the formal analysis of citation data cannot be considered as a one and only one basis for evaluation of the scientific impact [3]–[4]. However, the approach described in this article allows at least a rough comparison of the scientific impact of research institutions and to some extent even individual scientists working in different fields.

## References

[1] Science and Engineering Indicators 2004. National Science Foundation, May 04, 2004, Available on-line at: < http://www.nsf.gov/sbe/srs/seind04/ >  (accessed: October 26, 2004).

[2] ISIHighlyCited.Com, < http://isihighlycited.com/ > (accessed: October 26, 2004).

[3] Garfield, E.: Citation measures used as an objective estimate of creativity. Current Contents, #26, August 26, 1970 (see also: Garfield, E.: Essays of an Information Scientist, vol. 1, pp.120-121, 1962-73)

[4] Garfield, E.: Citation frequency as a measure of research activity and performance. Current Contents, #5, January 31, 1973 (see also: Garfield, E.: Essays of an Information Scientist, vol. 1, pp. 406-408, 1962-73)

[5] Garfield, E.: Is the ratio between number of citations and publications cited a true constant? Current Contents, #6, February 9, 1976 (see also: Garfield, E.: Essays of an Information Scientist, vol. 2, pp. 419-425, 1974-76)

[6] Garfield, E., and Welljarns-Dorof, A.: Citation data: their use as quantitative indicators for science and technology evaluation and policy-making. Science and Public Policy 19(5):32 1-7, October 1992

[7] Most-Cited Researchers 1983-2002, < http://www.sciencewatch.com/sept-oct2003/sw_sept-oct2003_page2.htm > (accessed: December 8, 2004).

Table 1. Comparison of the numbers of citations in different fields of science. Based on the data from *Science and Engineering Indicators 2004. National Science Foundation, May 04, 2004.*

| Field | Average ratio of citation number to the number of citations in mathematics | 1992 | | 1994 | | 1996 | | 1997 | | 1999 | | 2001 | |
|---|---|---|---|---|---|---|---|---|---|---|---|---|---|
| | | number of citations | ratio to maths | number of citations | ratio to maths | number of citations | ratio to maths | number of citations | ratio to maths | number of citations | ratio to maths | number of citations | ratio to maths |
| Clinical medicine | **78** | 475793 | 69 | 516665 | 78 | 554332 | 80 | 574859 | 90 | 584330 | 78 | 589762 | 76 |
| Biomedical research | **78** | 460148 | 67 | 518304 | 78 | 562361 | 81 | 572122 | 89 | 594596 | 79 | 568328 | 73 |
| Biology | **8** | 52535 | 8 | 57825 | 9 | 58649 | 8 | 58130 | 9 | 56981 | 8 | 57899 | 7 |
| Chemistry | **15** | 88010 | 13 | 96827 | 15 | 105960 | 15 | 105762 | 16 | 110927 | 15 | 109703 | 14 |
| Physics | **19** | 137922 | 20 | 141653 | 21 | 138417 | 20 | 131958 | 21 | 125968 | 17 | 120593 | 15 |
| Earth/space sciences | **9** | 55086 | 5 | 58818 | 9 | 71230 | 10 | 73507 | 11 | 83053 | 11 | 82614 | 11 |
| Engineering/technology | **5** | 32680 | 5 | 35189 | 5 | 33664 | 5 | 32958 | 5 | 34001 | 5 | 36809 | 5 |
| Mathematics | **1** | 6858 | 1 | 6631 | 1 | 6961 | 1 | 6418 | 1 | 7520 | 1 | 7794 | 1 |
| Social/behavioral sciences | **13** | 80282 | 12 | 84353 | 13 | 93032 | 13 | 93187 | 15 | 99481 | 13 | 104793 | 13 |